\documentclass{article}

\usepackage[utf8]{inputenc} 
\usepackage[T1]{fontenc} 
\usepackage{lmodern} 

\usepackage[top=3cm, includefoot]{geometry}

\usepackage{amsmath,amssymb,amsfonts,amsthm}
\usepackage{accents}
\usepackage{dsfont}

\usepackage{xifthen} 
\usepackage{graphicx} 
\usepackage{xcolor} 
 \usepackage{multirow}

\usepackage{caption}
\usepackage{subcaption}
\usepackage{tabularx}

\usepackage{tikz}
\usetikzlibrary{calc}
\usetikzlibrary{arrows.meta}
\usetikzlibrary{lindenmayersystems}
\usetikzlibrary{backgrounds}
\usepackage{array}
\usepackage{comment}
\usetikzlibrary{decorations.text, decorations.pathreplacing, shapes.geometric, arrows.meta, calc, shadows}

\usepackage{standalone}

\usepackage{pgfplots}
\pgfplotsset{compat=1.18}

\usepackage{pythonhighlight}

\newtheorem{theorem}{Theorem}[section]

\newtheorem{remark}[theorem]{Remark}
\newtheorem{proposition}[theorem]{Proposition}
\newtheorem{conjecture}[theorem]{Conjecture}

\usepackage{hyperref}
\usepackage{cleveref}

\usepackage{pdfpages}

\usepackage{todonotes}

\usepackage[all]{xy} 
\newcommand{\RR}{\mathbb{R}}
\newcommand{\NN}{\mathbb{N}}
\newcommand{\ZZ}{\mathbb{Z}}

\newcommand{\ee}{\mathrm{e}}
\newcommand{\OO}{\mathcal{O}}

\newcommand{\gee}{>}
\newcommand{\lee}{<}

\providecommand{\Norm}[2][]{\left\lVert#2\right\rVert\ifthenelse{\isempty{#1}}{}{_{#1}}}
\providecommand{\norm}[2][]{\lVert#2\rVert\ifthenelse{\isempty{#1}}{}{_{#1}}}
\providecommand{\bignorm}[2][]{\bigl\lVert#2\bigr\rVert\ifthenelse{\isempty{#1}}{}{_{#1}}}
\providecommand{\Bignorm}[2][]{\Bigl\lVert#2\Bigr\rVert\ifthenelse{\isempty{#1}}{}{_{#1}}}
\providecommand{\biggnorm}[2][]{\biggl\lVert#2\biggr\rVert\ifthenelse{\isempty{#1}}{}{_{#1}}}
\providecommand{\Biggnorm}[2][]{\Biggl\lVert#2\Biggr\rVert\ifthenelse{\isempty{#1}}{}{_{#1}}}

\makeatletter
\def\moverlay{\mathpalette\mov@rlay}
\def\mov@rlay#1#2{\leavevmode\vtop{%
   \baselineskip\z@skip \lineskiplimit-\maxdimen
   \ialign{\hfil$\m@th#1##$\hfil\cr#2\crcr}}}
\newcommand{\charfusion}[3][\mathord]{
    #1{\ifx#1\mathop\vphantom{#2}\fi
        \mathpalette\mov@rlay{#2\cr#3}
      }
    \ifx#1\mathop\expandafter\displaylimits\fi}
\makeatother

\newcommand{\bigcupdot}{\charfusion[\mathop]{\bigcup}{\cdot}}

\counterwithin{figure}{section}

\usepackage{graphicx}

\title{A family of non Minkowski measurable fractals in $\RR^2$}
\author{Uta Freiberg, Jonas Lippold}
\date{\today}

\pgfdeclarelindenmayersystem{Cantor dust}{
			\symbol{S}{\pgfpathrectangle{\pgfpointorigin}{\pgfpoint{1.0/1.0*\pgflsystemcurrentstep}{1.0/1.0*\pgflsystemcurrentstep}}}
			\symbol{U}{\pgftransformyshift{1.0/1.0*\pgflsystemcurrentstep}}
			\symbol{R}{\pgftransformxshift{1.0/1.0*\pgflsystemcurrentstep}}
			\rule{S -> [SRRS][UUSRRS]}
			\rule{U -> UUU}
			\rule{R -> RRR}
		}

\pgfdeclarelindenmayersystem{Cantor dust Epsilon}{
			\symbol{S}{\filldraw[color=black!20, fill=black!20, very thin, rounded corners=\myeps*\mystepsize](-\myeps,-\myeps) rectangle (\myepsend,\myepsend);}
			\symbol{U}{\pgftransformyshift{1.0/1.0*\pgflsystemcurrentstep}}
			\symbol{R}{\pgftransformxshift{1.0/1.0*\pgflsystemcurrentstep}}
			\rule{S -> [SRRS][UUSRRS]}
			\rule{U -> UUU}
			\rule{R -> RRR}
		}

\begin{document}

\maketitle

\begin{abstract}
    A long-standing conjecture of Lapidus asserts that, under certain conditions, a self-similar fractal set is not Minkowski measurable if and only if it is of lattice-type. For self-similar sets in $\RR$, the Lapidus conjecture has been confirmed. However, in higher dimensions, it remains unclear whether all lattice-type self-similar sets are not Minkowski measurable. This work presents a family of lattice-type subsets in $\RR^2$ that are not Minkowski measurable, hence providing further support for the conjecture. Furthermore, an argument is presented to illustrate why these sets are not covered by previous results.
\end{abstract}

\section{Introduction}
The Minkowski content is a useful tool for characterizing fractal sets beyond their (Hausdorff or Minkowski) dimension, particularly for distinguishing between sets of the same dimension. Therefore, it is of interest to determine which sets are Minkowski measurable and which are not.

Significant progress has been made in understanding Minkowski measurability, especially for self-similar sets generated by an iterated function system (IFS) satisfying the open set condition (OSC). In the 1990s, Lapidus conjectured that, in this setting, Minkowski measurability is generally equivalent to the IFS being non-lattice, an algebraic property that can be easily verified from the scaling ratios of the similarities in the IFS.

While the case of non-lattice sets being Minkowski measurable has been settled in $\RR^d$ by Gatzouras \cite{Gatzouras}, the lattice case remains more challenging. For $d=1$, the Lapidus conjecture has been fully formulated and proven, see \cite{Falconer:Minkowski, ReellerFall:neu}. However, for $d\geq 2$, the question remains open and is an active area of research. For a detailed survey, see \cite{Kombrink:survay:on:MM}.

Under certain technical assumptions, a sufficient condition for non-Minkowski measurability in higher dimensions is given by the pluriphase condition, see  \cite{Pluriphase}. This, along with the aforementioned results, is based on the renewal theorem. A different approach is the direct analysis of exact tube formulas, focusing on complex dimensions and fractal zeta functions. Recent research has shifted toward this powerful tools, but it often leads to highly involved calculations. 

This paper investigates a family of lattice self-similar Cantor dusts in $\RR^2$ and demonstrates, using elementary geometric arguments, that these sets are not Minkowski measurable. Additionally, an argument is presented suggesting why these sets fail to satisfy the pluriphase condition and have thus not been covered by previous results given in \cite{Pluriphase}. More specifically, the family depends on a real parameter $r > 2$. For $r\geq 30$, non-Minkowski measurability is  proven rigorously, while for $r \in (2,30)$, the question is reduced to an inequality that is verified numerically. 

The paper is organized as follows. In Section \ref{sec:preliminaries} preliminaries are introduced (see \ref{ssec: M.-mb.} - \ref{ssec: Pl_cond} ) and a summary of key results from the literature is provided (see \ref{ssec: known_results}). The family of Cantor dusts $\{C^r:r>2\}$ is introduced in Section \ref{sec: C^r}  and its relation to the pluriphase condition is discussed in Subsection \ref{ssec: pluri_for_c^r}. The main result, addressing the case $r\geq 30$, is stated and proved in Subsection \ref{ssec: case_r>30}, while the case $r \in (2,30)$ is analysed numerically in Subsection \ref{ssec: case_r<30}.

\section{Preliminaries}
\label{sec:preliminaries}

Essential terminology is introduced in order to state and prove the main theorems. For a detailed discussion of the basics of fractal geometry, the book \emph{Fractal Geometry}  by K. Falconer is recommended, see \cite{Falconer:FracGeo}.

\subsection{Minkowski measurability}
\label{ssec: M.-mb.}

Let  $(\RR^d, \| \cdot \|)$ be the $d$-dimensional Euclidean space and denote by 
$$\mathcal{H}(\RR^d) := \{ K \subset \RR^d \mid K \not= \emptyset \ \text{and} \ K \ \text{is compact} \}$$
the corresponding \emph{Hausdorff space}. For $A \in \mathcal{H}(\RR^d)$,  $\varepsilon > 0$ define the \emph{$\varepsilon$- parallel set} of $A$ as $$ A_\varepsilon := \{x \in \RR^d \mid \min _{a\in A} \| x-a\| \leq \varepsilon \}.$$
 Denote by $\lambda^d$ the $d$-dimensional Lebesgue measure. If the \emph{Minkowski dimension} $$\dim_{\mathcal{M}}(A):=d-\lim_{\varepsilon \searrow 0} \frac{\ln \lambda^d(A_\varepsilon)}{\ln\varepsilon}$$ exists, then the \emph{Minkowski content} of $A$ is given by
 $$\mathcal{M}(A):= \lim_{\varepsilon \searrow 0} \frac{\lambda^d (A_\varepsilon)}{\varepsilon^{d-\dim_{\mathcal{M}}(A)}}.$$

The set $A$ is called $\emph{Minkowski measurable}$ if $\mathcal{M}(A)$ exists and is both positive and finite.

\subsection{Self-similar sets, open set condition, and the (non-)lattice case}

Let $S=\{ S_1, \ldots S_N\}$, with $N\geq 2$, denote an iterated function system (IFS) consisting of contractive similarities $S_1, \ldots S_N$ acting on  $\RR^d$. Such a system is called a \emph{self-similar system}. Define the corresponding (set-valued) map $ \mathbf{S} : \mathcal{H}(\RR^d) \to \mathcal{H}(\RR^d) $ by
\begin{equation}
\label{eq:fett_S}
    \mathbf{S}A:=\bigcup_{i=1}^{N} S_i(A), \quad A \in \mathcal{H}(\RR^d).
\end{equation}

It is well known (\cite{Hut}) that $ \mathbf{S} : \mathcal{H}(\RR^d) \to \mathcal{H}(\RR^d) $ has a unique fixed point $F$ which is called the \emph{self-similar attractor} of the self-similar system $S$. A set $A \in \mathcal{H}(\RR^d)$ is called a \emph{self-similar set} if it is the attractor of some self-similar system. 

The IFS $S$ satisfies the \emph{open set condition} (OSC) if there exits a non-empty bounded open set $\mathcal{O} \subset \RR^d $ such that  
\begin{equation}
\label{eq:osc}
    \bigcupdot_{i=1}^N S_i(\mathcal{O}) \subseteq \mathcal{O},
\end{equation}
where $\bigcupdot$ denotes the pairwise disjoint union. Any non-empty open set $\OO$ satisfying \eqref{eq:osc} is called a \emph{feasible open set} for $S$. If, in addition, $\OO \cap F \not = \emptyset $, then $\OO $ is called a \emph{strong feasible open set}. Note that it was shown in \cite{Schief} that if a self-similar system satisfies (OSC), then it possesses a strong feasible open set. 

Denote the scaling ratio of $S_i$ by $r_i$, $i=1,\ldots N$. The IFS $S$ is said to be \emph{lattice} if $\{ \ln r_1, \ldots , \ln r_N \}$ generates a discrete subgroup of $(\RR, +).$ Otherwise, $S$ is said to be  \emph{non-lattice}. If $S$ is lattice, then there exists a largest $a>0$ such that $\{ \ln r_1, \ldots , \ln r_N \} \subseteq a\ZZ$ and $\ee^a$ is called the \emph{base} of $S$.

\subsection{Non-triviality and projection condition}
Let $F \in \mathcal{H}(\RR^d)$ be the attractor of a self-similar system $S =\{S_1, \ldots S_N \}$ satisfying (OSC). $F$ is called \emph{non-trivial} if there exists a feasible open set $\OO$ such that
\begin{equation}
    \OO \not \subseteq \overline{\mathbf{S}\OO},
\end{equation}
where $\overline{B}$ denotes the closure of $B \subseteq \RR^d$; otherwise, $F$ is called \emph{trivial}.

$F$ is non-trivial if and only if $F$ has an empty interior, which is equivalent to $F$ having a Minkowski dimension strictly less than $d$, as shown in \cite{Pears:tilings}. In particular, non-triviality is independent of the choice of the feasible set $\OO$.

Let $\pi_A$ denote the \emph{metric projection} onto $A$ for $A \in \mathcal{H}(\RR^d)$, which is defined for points $x \in \RR^d$ that have a unique nearest point $y$ in $A$ by $$\pi_F(x):= y.$$ The set $\OO$ is said to satisfy the \emph{projection condition} if $$S_i\OO \subseteq \overline{\pi_F^{-1}(S_i F)}$$ for all $i=1,\ldots N$. It is worth noting that, as long as (OSC) holds, it is always possible to find a strong feasible open set that satisfies the projection condition (see \cite{Winter:MinContent}).

\subsection{Pluriphase condition}
\label{ssec: Pl_cond}
For a given IFS with non-trivial attractor $F \in \mathcal{H}(\RR^d)$ and a fixed feasible open set $\OO$, define $$\Gamma := \Gamma_{\OO} := \OO \setminus \mathbf{S}\OO$$ and $$g:=\sup\limits_{x\in\Gamma} \{\min\limits_{y\in F}\{\|x-y\|\}\}.$$ The set $F$ is said to be \emph{pluriphase  with respect to} $\Gamma_{\OO}$ if there exists a finite partition $0=: a_0 < a_1<\ldots <a_M:= g < \infty$ of the interval $(0,\infty)$ and constants $\kappa_{m,k} \in \RR, \ m=1,\ldots ,M, \ k =0,\ldots ,d$, such that for all $\varepsilon>0$
\begin{equation}
    \lambda^{d} (F_\varepsilon \cap \Gamma)=\sum_{m=1}^M \mathds{1}_{(a_{m-1},a_m]}(\varepsilon) \sum_{k=0}^d \kappa_{m,k}\varepsilon^{d-k} + \mathds{1}_{(g,\infty)}(\varepsilon)\lambda^{d} (\Gamma).
\end{equation}
Hereby $\mathds{1}_A$ denotes the indicator function of a set $A$.

\subsection{Known results on Minkowski measurability of self-similar sets}
\label{ssec: known_results}
\begin{figure}

    \centering
    \begin{tikzpicture}[%
    block/.style={
      draw,
      rectangle,
      rounded corners=5pt,
      fill=orange!20,
      minimum width=2.5cm,
      minimum height=4cm,
      align=center,
      drop shadow={shadow xshift=1pt,shadow yshift=-1pt,opacity=0.4}
    },
    block2/.style={
      draw,
      rectangle,
      rounded corners=5pt,
      fill=teal!20,
      minimum width=2.5cm,
      minimum height=4cm,
      align=center,
      drop shadow={shadow xshift=1pt,shadow yshift=-1pt,opacity=0.4}
    },
    arrow/.style={
      double distance=2pt,
      -{Latex[length=3mm]},
      thick,
      draw=black!50!green,
      shorten >=3pt,
      shorten <=3pt
    },
     arrow2/.style={
      double distance=2pt,
      -{Latex[length=3mm]},
      thick,
      draw=black!20!blue,
      shorten >=3pt,
      shorten <=3pt
    }
  ]

  \node[block]  (A) at (-6,0) {non-\\lattice};
  \node[block] (B) at (6,0)  {Minkowski\\measurable};

  \draw[arrow]
    ([yshift=45pt]A.east) -- ([yshift=45pt]B.west)
    node[midway, above, font=\small\itshape] {$\mathbb{R}^{\ge1}$}
    node[midway, below, font=\small\itshape] {\emph{Gatzouras 2000}};

  \draw[arrow2]
    ([yshift=0pt]B.west) -- ([yshift=0pt]A.east)
    node[midway, above, font=\small\itshape] {$\mathbb{R}^1,\ \dim_{\mathcal M}<1$}
    node[midway, below, font=\small\itshape] {\emph{Falconer 1995, Kombrink and Winter 2020}};

  \draw[arrow2]
    ([yshift=-45pt]B.west) -- ([yshift=-45pt]A.east)
    node[midway, above, font=\small\itshape] {$\mathbb{R}^{\ge2},\,\dim_{\mathcal M}\notin\mathbb N$, \emph{pluriphase\,w.r.t.}\,\(\mathcal O\)}
    node[midway, below, font=\small\itshape] {\emph{Kombrink, Pearse and Winter 2016}};

\end{tikzpicture}
    \vspace{0.3cm}
    \caption{Overview of the relationships between non-lattice and Minkowski measurability of self-similar sets satisfying the (OSC).} 
    \label{fig:Schaubild}
\end{figure}
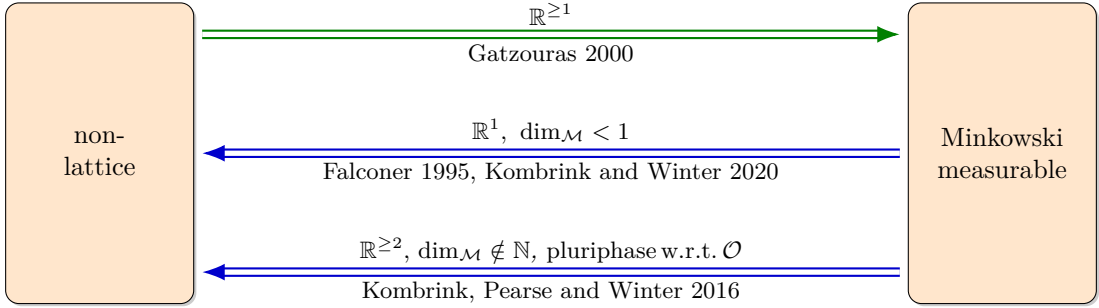

For self-similar sets in $\RR^d$, the most important results on Minkowski measurability are presented in order to understand the connection between Minkowski measurability and the non-lattice property. See also Figure \ref{fig:Schaubild} for a schematic overview. 
\begin{theorem}[Gatzouras, \cite{Gatzouras}] \label{thm: nonlattice_folgt_M-mb}
Let be $F \in \mathcal{H}(\RR^d)$ the self-similar attractor of an non-lattice, (OSC) self-similar system. Then $F$ is  Minkowski measurable.
\end{theorem}

\begin{theorem}
    [Falconer, \cite{Falconer:Minkowski} and Kombrink, Winter, \cite{ReellerFall:neu}] \label{thm: R^1_lattice_folgt_nicht-M-mb} Let be $F \in \mathcal{H}(\RR)$ the self similar attractor of an lattice, (OSC) self-similar system with $\dim_{\mathcal{M}}F \lee 1$. Then $F$ is not Minkowski measurable.
\end{theorem}

Theorems \ref{thm: nonlattice_folgt_M-mb} and \ref{thm: R^1_lattice_folgt_nicht-M-mb} together give a characterization of Minkowski measurability for self-similar sets generated by an (OSC) self-similar system $S$ in $\RR$. If these sets have Minkowski dimension less than one, then they are Minkowski measurable if and only if $S$ is a non-lattice IFS. It remains unclear whether this equivalence holds in higher dimensions. Theorem \ref{thm: pluriphase-case} provides a first step in this direction.

\begin{theorem}[Kombrink, Pearse, Winter, \cite{Pluriphase}] \label{thm: pluriphase-case}
Let be $F \in \mathcal{H}(\RR^d)$ with $\dim_{\mathcal{M}} F \not \in \NN$ the self-similar attractor of an non-lattice, (OSC) self-similar system. Suppose that there exists a strong feasible open set $\mathcal{O}$ satisfying the projection condition such that $F$ is pluriphase with respect to $\Gamma_{\OO}$. Then $F$ is not Minkowski measurable.
\end{theorem}

Note that it is essential to exclude sets with an integer Minkowski dimension from the Lapidus conjecture.

\section{The Cantor dust $C^r$}

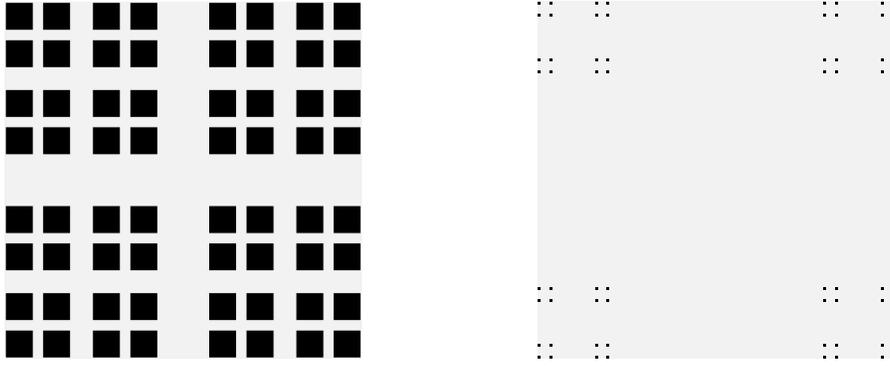
\begin{figure}
    \centering
    \begin{tabular}{ c @{\qquad} c @{\qquad} c }

\begin{tikzpicture}[scale=0.0138]
		\filldraw[color=gray!10, fill=black!100] (0,0) -- (0,343) -- (343,343) -- (343,0) -- (0,0);
		
		\filldraw[color=gray!10, fill=gray!10] (3*49,0) -- (4*49,0) -- (4*49,343) -- (3*49,343)--(3*49,0);
		\filldraw[color=gray!10, fill=gray!10] (9*7,0) -- (12*7,0) -- (12*7,343) -- (9*7,343)--(9*7,0);
		\filldraw[color=gray!10, fill=gray!10] (4*49+9*7,0) -- (4*49+12*7,0) -- (4*49+12*7,343) -- (4*49+9*7,343)--(4*49+9*7,0);
		\filldraw[color=gray!10, fill=gray!10] (9*3,0) -- (9*4,0) -- (9*4,343) -- (9*3,343)--(9*3,0);	
		\filldraw[color=gray!10, fill=gray!10] (4*49+9*3,0) -- (4*49+9*4,0) -- (4*49+9*4,343) -- (4*49+9*3,343)--(4*49+9*3,0);
		\filldraw[color=gray!10, fill=gray!10] (12*7+9*3,0) -- (12*7+9*4,0) -- (12*7+9*4,343) -- (12*7+9*3,343)--(12*7+9*3,0);
		\filldraw[color=gray!10, fill=gray!10] (4*49+12*7+9*3,0) -- (4*49+12*7+9*4,0) -- (4*49+12*7+9*4,343) -- (4*49+12*7+9*3,343)--(4*49+12*7+9*3,0);

		\filldraw[color=gray!10, fill=gray!10] (0,3*49) -- (343,3*49) -- (343,4*49) -- (0,4*49) -- (0,3*49);
		\filldraw[color=gray!10, fill=gray!10] (0,9*7) -- (0,12*7) -- (343,12*7) -- (343,9*7)--(0,9*7);
		\filldraw[color=gray!10, fill=gray!10] (0,4*49+9*7) -- (0,4*49+12*7) -- (343,4*49+12*7) -- (343,4*49+9*7)--(0,4*49+9*7);
		\filldraw[color=gray!10, fill=gray!10] (0,9*3) -- (0,9*4) -- (343,9*4) -- (343,9*3)--(0,9*3);	
		\filldraw[color=gray!10, fill=gray!10] (0,4*49+9*3) -- (0,4*49+9*4) -- (343,4*49+9*4) -- (343,4*49+9*3)--(0,4*49+9*3);
		\filldraw[color=gray!10, fill=gray!10] (0,12*7+9*3) -- (0,12*7+9*4) -- (343,12*7+9*4) -- (343,12*7+9*3)--(0,12*7+9*3);
		\filldraw[color=gray!10, fill=gray!10] (0,4*49+12*7+9*3) -- (0,4*49+12*7+9*4) -- (343,4*49+12*7+9*4) -- (343,4*49+12*7+9*3)--(0,4*49+12*7+9*3);	
	\end{tikzpicture}
 &&
 \pgfdeclarelindenmayersystem{Cantor dust 5}{
    \symbol{S}{\pgfpathrectangle{\pgfpointorigin}
        {\pgfpoint{1.0/1.0*\pgflsystemcurrentstep}{1.0/1.0*\pgflsystemcurrentstep}} \pgfusepath{fill}}
    \symbol{U}{\pgftransformyshift{1.0/1.0*\pgflsystemcurrentstep}}
    \symbol{R}{\pgftransformxshift{1.0/1.0*\pgflsystemcurrentstep}}
    \rule{S -> [SRRRRS][UUUUSRRRRS]}
    \rule{U -> UUUUU}
    \rule{R -> RRRRR}
}

\begin{tikzpicture}[scale=0.5]
    \begin{scope}[on background layer]
        \fill[gray!10] (0,0) rectangle (9.5,9.5);
    \end{scope}
    
    \fill[lindenmayer system={Cantor dust 5, axiom=S, step=2.16pt, order=3}] lindenmayer system;
\end{tikzpicture}

\end{tabular}
    \caption{$C^r$ for $r=7/3$ (left) and $r=5$ (right).}
    \label{fig:C^r}
\end{figure}

\label{sec: C^r}

Let $r\in (2,\infty)$ be a parameter and consider the maps $ S^r_1,S^r_2,S^r_3,S^r_4 : \RR^2 \to \RR^2$ given by 

$$\ S^r_1(x)=\frac{1}{r}x, \hspace{3cm} \ S^r_2(x)=\frac{1}{r}x+\begin{pmatrix} \frac{r-1}{r} \\ 0 \end{pmatrix},$$ $$
 \ S^r_3(x)=\frac{1}{r}x+\begin{pmatrix} \frac{r-1}{r} \\ \frac{r-1}{r} \end{pmatrix}, \hspace{1.419cm} \ S^r_4(x)=\frac{1}{r}x+\begin{pmatrix} 0 \\ \frac{r-1}{r} \end{pmatrix}.$$ 
Denote the corresponding IFS by
\begin{equation}
\label{eq:IFS_S^r}
    S^r :=\{ S^r_1,S^r_2,S^r_3,S^r_4\},
\end{equation} as well as the associated self-similar attractor by  $C^r$, see Figure \ref{fig:C^r}.
It is clear that $S^r$ satisfies (OSC), with $ \mathcal{O} = (0,1)\times (0,1)$ as a feasible set. Moreover, $S^r$ is lattice with base $r$ and $C^r$ has Minkowski (and Hausdorff) dimension $D_r=\frac{\ln 4}{\ln r}$; see \cite{Hut}.

Let $W:=\{1,2,3,4\}$ be an alphabet. For a word $w=w_1w_2\dots w_n \in W^n$ of lengths $n$, define $S_w : \RR^2 \to \RR^2 $ by $$S_w:=S_{w_1} \circ \ldots \circ S_{w_n}.$$  Then the  $n$-th construction step of $C^r$ is given by $$\bigcup_{w\in W^n}S_w([0,1]\times[0,1])$$ and consists of $4^n$ disjoint squares, each of side length $r^{-1}$.

\begin{figure}
     \centering
	\begin{tabular}{ c @{\qquad \qquad \qquad \qquad} c }
         \includegraphics[width=0.30\linewidth]{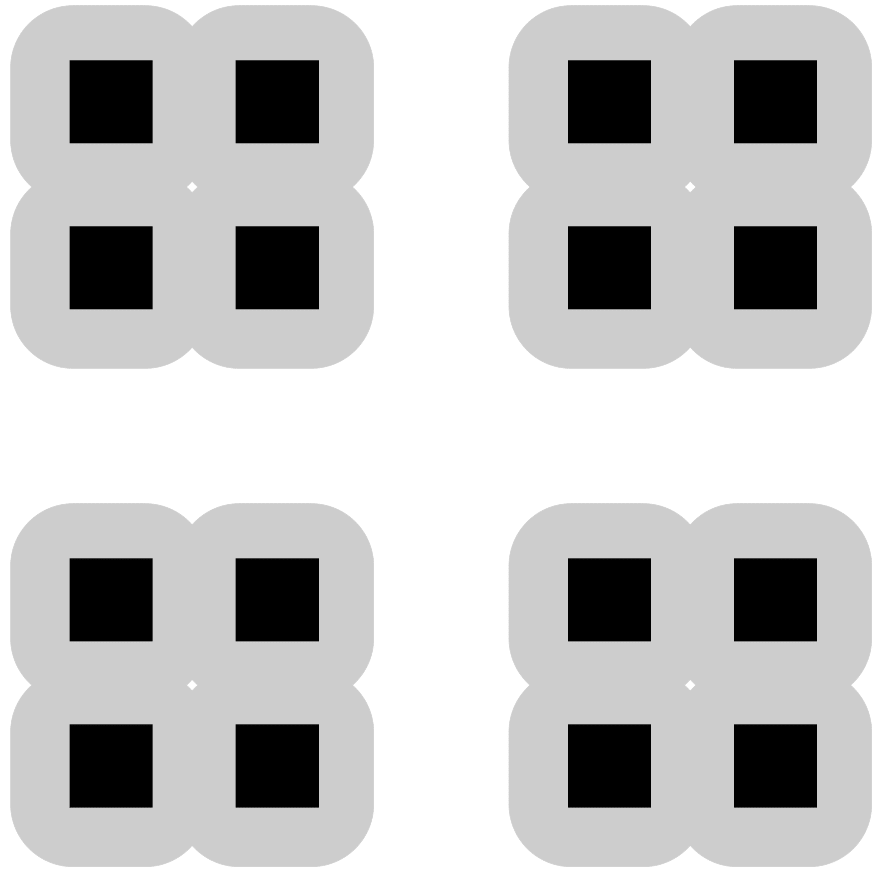}
        &
        \includegraphics[width=0.30\linewidth]{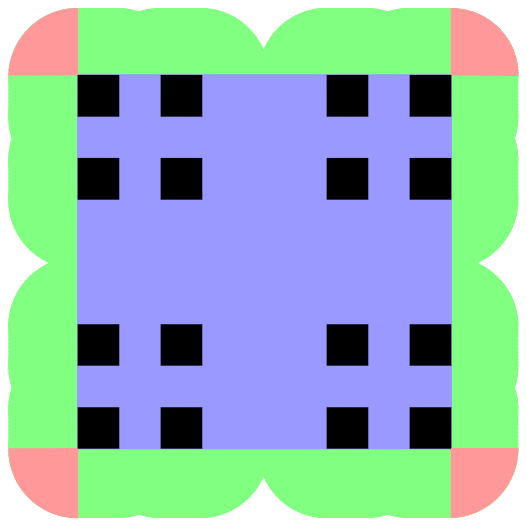}
    \\
		\phantom{x}& 
		\\
		\small (a)& 
		\small (b)
	\end{tabular}
    
     \caption{(a) An $\varepsilon$- neighbourhood of $C^r$, where each component contains no holes. (b) A single component of the $\varepsilon$-neighbourhood of $C^r$ for $\varepsilon$ satisfying Condition (\ref{eqn:geeignete_Folgenwahl}).}
     \label{fig:BGR_neighbourhood} 
    
 \end{figure}

The $\varepsilon$-neighbourhood of $C^r$, where $\varepsilon$ satisfies
\begin{equation}
\label{eqn:geeignete_Folgenwahl}
\sqrt{\frac{1}{2}}\cdot\frac{r-2}{r}r^{-n}\leq \varepsilon \leq \frac{r-2}{2}r^{-n},
\end{equation}
consists of $4^n$ congruent components, where the intersection of two different components is at most one dimensional and each component contains no holes. This is illustrated in Figure \ref{fig:BGR_neighbourhood}, where each component can be divided into three parts:

\begin{itemize}
    \item $4$ quarter circles of radius $\varepsilon$, denoted by $R(\varepsilon)$, marked in red,
    \item one square of side length $r^{-1}$, denoted by $B(\varepsilon)$, marked in blue,
    \item $4$ parts with fractal-like boundary  $G(\varepsilon)$, marked in green. 
    \end{itemize}  

This decomposition provides a structured approach for handling the area $\lambda^2(C^r_\varepsilon)$ in the upcoming calculations on the Minkowski measurability of $C^r$. Note that, for simplicity, the same notation is used for different but congruent figures.
 
 \subsection{The pluriphase condition for $C^r$}
 \label{ssec: pluri_for_c^r}
Before discussing the Minkowski measurability of $C^r$, this section demonstrates that the connection between the pluriphase condition and Minkowski measurability, as stated in Theorem \ref{thm: pluriphase-case}, may not be applicable. At least in the case of the canonical IFS $S^r$ given in (\ref{eq:IFS_S^r}) generating $C^r$, with $(0,1)\times(0,1)$ chosen as the feasible set, the attractor is not pluriphase. Denote by $\mathbf{S}^r$ the map corresponding to $S^r$, as defined in (\ref{eq:fett_S}).

\begin{proposition}
\label{prop:C^r_pluri}
Let $r\in (2,\infty)$. Then $C^r$ is not pluriphase w.r.t. $\Gamma = (0,1)\times(0,1)\setminus \mathbf{S}^r((0,1)\times(0,1))$.
\end{proposition}

\begin{proof}
Let $r>2$ and $\OO:=(0,1)\times(0,1)$. Then $\OO$ is a feasible set for $S^r$ for all $r>2$, and $\Gamma=\OO\setminus\mathbf{S}^r\OO$ forms a "cross", as depicted in Figure \ref{fig:Gamma}.

For $\varepsilon < \frac{r-2}{2r}$, denote by $P(\varepsilon)$ the area of $\Gamma \cap C^r_\varepsilon \setminus ([ r^{-1}, 1- r^{-1} ] \times [ r^{-1}, 1- r^{-1} ])$, as illustrated in Figure \ref{fig:Gamma}. Then, $$\lambda^2(C^r_{\varepsilon} \cap \Gamma) = P(\varepsilon) + \pi \varepsilon^2.$$

The function $P$ satisfies the recursion
\begin{equation*}
    P(r^{-1} \varepsilon) = 2 \cdot r^{-2}P(\varepsilon) + 8 \cdot r^{-2}\frac{\pi}{2} \varepsilon^2.
\end{equation*}
Thus, the following holds:

\begin{equation}
\label{eq:rekursionsformel}
    \begin{split}
\lambda^{2} (C^r_{\frac{1}{r}\varepsilon}\cap \Gamma) &=  2 r^{-2}P(\varepsilon) + 8 r^{-2}\frac{\pi}{2} \varepsilon^2 + r^{-2} \pi \varepsilon^2\\
&= 2 r^{-2} \Big[P(\varepsilon) + \frac{5}{2} \pi \varepsilon^2 \Big] \\
&= 2 r^{-2} \Big[\lambda^2(C^r_{\varepsilon} \cap \Gamma) + \frac{3}{2} \pi \varepsilon^2 \Big].
\end{split}
\end{equation}

 \begin{figure}
     \centering

	\begin{tabular}{ c @{\qquad} c @{\qquad} c }
         \includegraphics[width=0.25\linewidth]{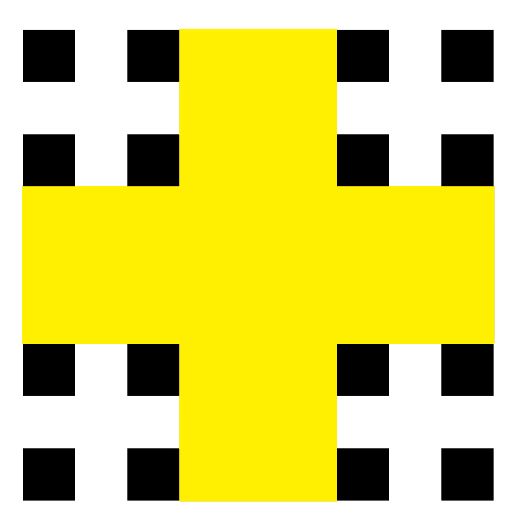}
        &
        \includegraphics[width=0.25\linewidth]{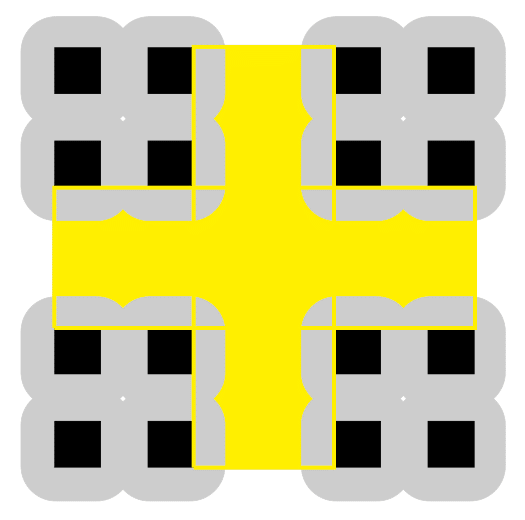}
        &
         \includegraphics[width=0.25\linewidth]{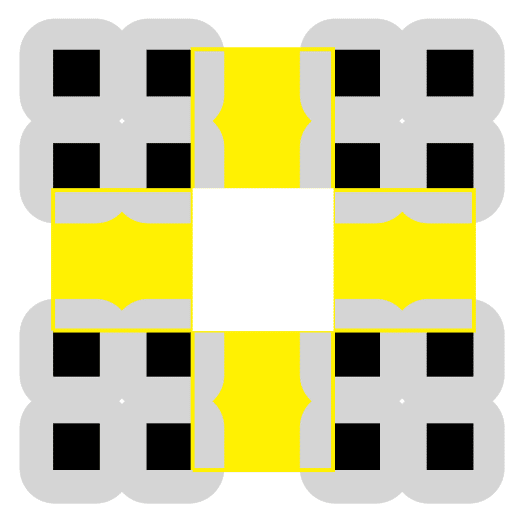}
    \\
		\phantom{x}& 
		\\
		\small $\Gamma$& 
		\small $\Gamma \cap C^r_\varepsilon$&
        \small $P(\varepsilon)$
	\end{tabular}
    
     \caption{Illustration of $\Gamma$ in yellow, along with the subset $\Gamma \cap C^r_\varepsilon$ and the function $P(\varepsilon)$ used in the proof of Proposition \ref{prop:C^r_pluri}.}
     \label{fig:Gamma}
 \end{figure}

Assume that $C^r$ is pluriphase w.r.t. $\Gamma$. Then there exists a $\varepsilon_0 > 0$ such that $q(\varepsilon):= \lambda^2(C^r_\varepsilon \cap \Gamma)$ is a polynomial of degree at most $2$ on $(0, \varepsilon_0]$. In other words, there exists constants $a,b,c \in \RR$ with $q(\varepsilon)= a\varepsilon^2 + b \varepsilon + c$. From Equation (\ref{eq:rekursionsformel}), for $\varepsilon < \varepsilon_0$ it follows that $$q(r^{-1}\varepsilon)=2r^{-2}\Big[q(\varepsilon)+\frac{3}{2}\pi \varepsilon^2\Big],$$ and hence
$$a\Big(\frac{1}{r}\varepsilon \Big)^2 + b \Big(\frac{1}{r}\varepsilon\Big) + c = 2r^{-2} \Big(a\varepsilon^2 + b \varepsilon + c + \frac{3}{2} \pi \varepsilon^2 \Big).$$

This yields $a=-3\pi,$ $b=c=0$. Consequently, $q(\varepsilon) = -3\pi \varepsilon^2$, which contradicts $q(\varepsilon)=\lambda^2(C^r_\varepsilon \cap \Gamma)>0$. 

\end{proof}

\begin{remark} There exist various (OSC) IFSs for which $C^r$ is the attractor. By definition, the pluriphase condition depends on both the choice of the IFS and the feasible set. However, at least for the canonical choice of an IFS generating $C^r$ and the feasible set $(0,1)\times(0,1)$, Theorem \ref{thm: pluriphase-case} does not provide a straightforward way to prove the non-Minkowski measurability of $C^r$.
\end{remark}

 \section{On the Minkowski measurability of $C^r$ }

In view of the Lapidus-conjecture, the aim is to prove that $C^r$ is not Minkowski measurable for all $r\gee 2$; that is, to show that the limit $$\lim_{\varepsilon \searrow 0} \frac{\lambda^{2} (C^r_\varepsilon)}{\varepsilon^{2-D_r}}$$ does not exist. To establish this result, it will be beneficial to introduce, for two $r$-dependent null sequences $(\varepsilon_{1,n})_{n\in \NN}$ and $(\varepsilon_{2,n})_{n\in \NN}$, the notations $$R_{i,n}:=R(\varepsilon_{i,n}), \quad B_{i,n}:=B(\varepsilon_{i,n}), \quad G_{i,n}:=G(\varepsilon_{i,n}).$$

\subsection{The case when $r\geq 30$ (Main result)}
\label{ssec: case_r>30}

For sufficiently large $r$ (specifically, $r \geq 30$) and for $\varepsilon$ satisfying Equation \eqref{eqn:geeignete_Folgenwahl}, the area of the $\varepsilon$-parallel set of $C^r$ is dominated by the red parts $R(\varepsilon)$ (see Figure \ref{fig:BGR_neighbourhood}). Therefore, it is relatively easy to obtain good estimates for $\lambda^2(C_\varepsilon)$.

 \begin{theorem}
 \label{Thm:case_r_grater_30}
     Let $r \in [30, \infty)$. Then $C^r$ is not Minkowski measurable.
 \end{theorem}
 \begin{proof}

 Let $r\geq 30$ and consider two null sequences defined by $$\varepsilon_{1,n}:= \sqrt{\frac{1}{2}}r^{-n} \qquad \text{and} \qquad \varepsilon_{2,n}:=r^{-n}.$$ Since $r\geq30$ it follows that $\varepsilon_{1,n}$ and $\varepsilon_{2,n}$  satisfy Equation  \eqref{eqn:geeignete_Folgenwahl}. Hence, $\lambda^{2} (C^r_{\varepsilon_{1,n}})$ is computed as follows:
\begin{equation}
\label{eq:eps_1_paralellvolum}
     \lambda^{2} (C^r_{\varepsilon_{1,n}}) = \lambda^{2} (B_{1,n} \cup G_{1,n} \cup R_{1,n})= 4^n r^{-2n} + \lambda^{2} ( G_{1,n}) + 4^n\frac{\pi}{2} r^{-2n}.
\end{equation}

The next step is to verify the following lower bound for the green areas:
\begin{equation}
\label{eq:eps_1_lower_bound}
  \lambda^{2} (G_{1,n}) > 4^n\cdot 4\cdot\frac{3}{4}\cdot r^{-n}\cdot\sqrt{\frac{1}{2}} r^{-n}.  
\end{equation}
This is depicted in Figure \ref{fig:lower_bound_for_r_greater_30} (a) and (c).

\begin{figure}[th]
	\centering
	\begin{tabular}{ c @{\qquad} c @{\qquad} c }
		
	\begin{tikzpicture}[scale=0.4]
\filldraw[color=green!50, fill=green!40] (0,0) -- (9,0) -- (9,2) arc[start angle=0, end angle=40, radius= 1] coordinate(C) arc[start angle=-40, end angle=0, radius= 1] coordinate(D)-- (9,5) arc[start angle=0, end angle=25, radius= 6] coordinate(E) arc[start angle=-25, end angle=0, radius= 6] coordinate(F)--(9,12) arc[start angle=0, end angle=40, radius= 1] coordinate(G) arc[start angle=-40, end angle=0, radius= 1] coordinate(H)-- (9,15) -- (0,15);
\draw[loosely dotted] (0,0)--(9,15);
\draw[loosely dotted] (0,15)--(9,0);
	\end{tikzpicture}
	& 
	\begin{tikzpicture}[scale=0.4]
\filldraw[color=green!50, fill=green!40] (0,0) -- (9,0) arc[start angle=0, end angle=55.3, radius= 9] coordinate(C) arc[start angle=-55.3, end angle=0, radius= 9] coordinate(D)-- (9,15)--(0,15);
\draw[loosely dotted] (0,0)--(9,15);
\draw[loosely dotted] (0,15)--(9,0);
	\end{tikzpicture}
	&
	\begin{tikzpicture}[scale=0.4]
\filldraw[color=green!50, fill=green!40] (0,0) -- (9,0) --(0,15)--(0,0);
\filldraw[color=green!50, fill=green!40] (0,0) -- (9,15) --(0,15)--(0,0);
\draw[loosely dotted] (0,0)--(9,15);
\draw[loosely dotted] (0,15)--(9,0);
	\end{tikzpicture}
	
	\\
		\phantom{x}& &
		\\
		\small $G_{1,n}$& 
		\small $\tilde{G}_{1,n}$&
		\small $G'_{1,n}$
	\end{tabular}
	\caption{Representation of the lower bound $G'_{1,n}$ of $G_{1,n}$, taking the intermediate step $\tilde{G}_{1,n}$. The sets satisfy $G'_{1,n} \subseteq \tilde{G}_{1,n} \subseteq G_{1,n}$.}
	\label{fig:lower_bound_for_r_greater_30}
\end{figure}
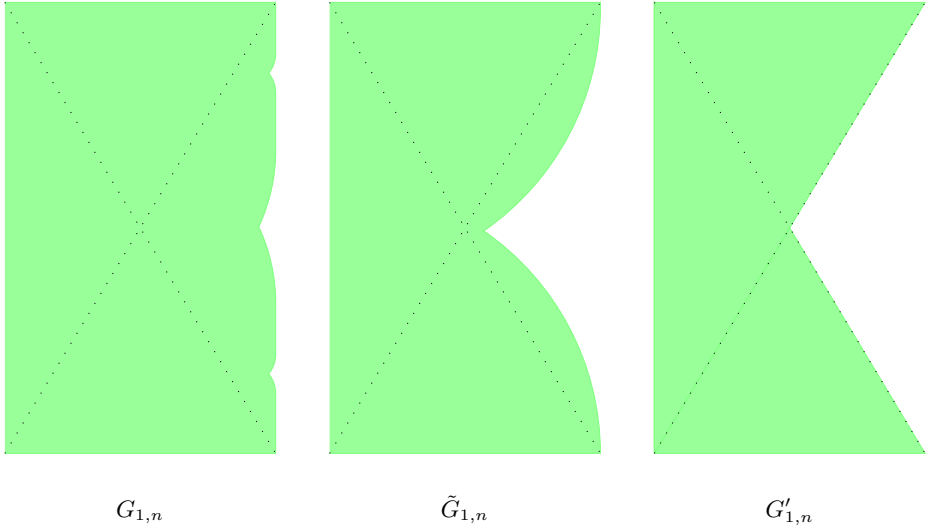
 
 To verify Equation \eqref{eq:eps_1_lower_bound}, denote by $V(B_{1,n})$ the vertices of $B_{1,n}$ and by $D(p, \varepsilon)$ the disk of radius $\varepsilon$ centered at p. Then, $$G_{1,n} = \bigcup_{p\in C^r} D(p, \varepsilon) \cap  G_{1,n} \subset \bigcup_{p \in V(B_{1,n+1})} D(p, \varepsilon) \cap G_{1,n} := \Tilde{G}_{1,n},$$ as depicted in Figure \ref{fig:lower_bound_for_r_greater_30}.

 Furthermore, note that $\lambda^2 (\Tilde{G}_{1,n})$ is greater than $4^n\cdot4\cdot\frac{3}{4}\cdot r^{-n} \cdot\sqrt{\frac{1}{2}} r^{-n}$ as shown in Figure \ref{fig:lower_bound_for_r_greater_30}. In this context, it should be observed that $2\cdot \varepsilon_{1,n}=2\cdot\sqrt{\frac{1}{2}} r^{-n} \geq r^{-n}$ and that for any rectangle $Q$ with side lengths $l_1$ and $l_2$ satisfying $l_1=2\cdot l_2$, one has
$\frac{3}{4}\lambda^{2} (Q)= \frac{3}{4}\cdot l_1\cdot l_2= \frac{3}{2} l_2^2 < \frac{\pi}{2} l_2^2 = \frac{1}{2} \lambda^{2} (D_{l_2}((0,0))) $, as illustrated in Figure \ref{fig:Worst_Case_Lower_bound_eps_1_1}.

\begin{figure}[th]
	\centering
	\begin{tabular}{ c @{\qquad} c @{\qquad} c }

	\begin{tikzpicture}[scale=3]

\filldraw[color=black!100, fill=black!10] (0,0) -- (2,0) --(2,1) arc[start angle=90, end angle=180, radius= 1] coordinate(C) arc[start angle=0, end angle=90, radius= 1] coordinate(D)-- (0,0) node[midway,left] {$l_1$};
\draw (0,0) -- (2,0) node[midway,below] {$l_2$};
\draw[loosely dotted] (0,0)--(2,1);
\draw[loosely dotted] (2,0)--(0,1);
	\end{tikzpicture}
	& &
	\begin{tikzpicture}[scale=3]
\filldraw[color=black!100, fill=black!10](0,0) -- (2,0) --(2,1) -- (1,0.5)--(0,1)--(0,0) node[midway,left] {$l_1$};
\draw (0,0) -- (2,0) node[midway,below] {$l_2$};
\draw[loosely dotted] (0,0)--(2,1);
\draw[loosely dotted] (2,0)--(0,1);
	\end{tikzpicture}
	
	\\
		\phantom{x}& &
		\\
		\small $\tilde{G}_{1,n}$& &
		\small $G'_{1,n}$
	\end{tabular}
    \caption{Even the smallest possible $\tilde{G}_{1,n}$ is still larger than the estimate  $G'_{1,n}$ of $G_{1,n}$.}
	\label{fig:Worst_Case_Lower_bound_eps_1_1}
\end{figure}
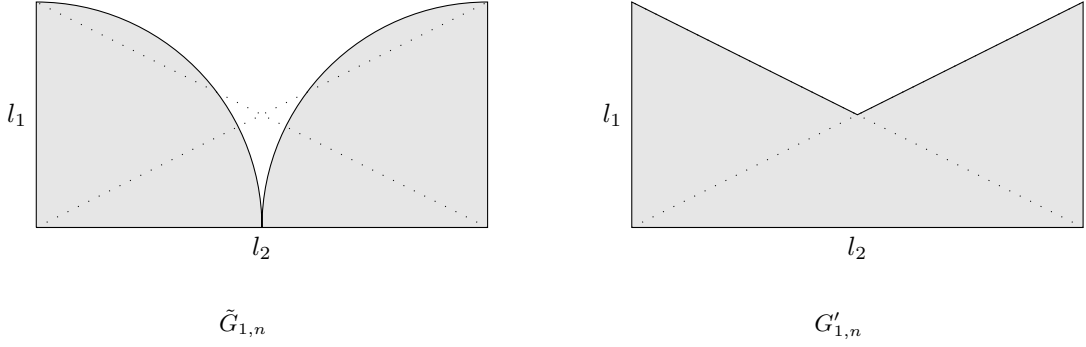

Using the equations (\ref{eq:eps_1_paralellvolum}) and (\ref{eq:eps_1_lower_bound})  it follows for all $n \in \NN $ 
\begin{equation}
\label{eqn:H1_vom_ersten_Satz}
\begin{split}
\frac{\lambda^{2} (C^r_{\varepsilon_{1,n}})}{\varepsilon_{1,n}^{2-D_r}} &>  \frac{4^n r^{-2n}\Big(1 + 3\sqrt{\frac{1}{2}} + \frac{\pi}{2}\Big)}{\Big(\sqrt{\frac{1}{2}} r^{-n}\Big)^{2}\Big(\sqrt{\frac{1}{2}} r^{-n}\Big)^{-\log_r(4)}} \\\ \\\
&= \frac{1 + 3\sqrt{\frac{1}{2}} + \frac{\pi}{2}}{\frac{1}{2}\sqrt{2}^{\log_r(4)}}.
\end{split}
\end{equation}

On the other hand, to estimate $\lambda^{2} (C^r_{\varepsilon_{2,n}})=\lambda^{2} (B_{2,n} \cup G_{2,n} \cup R_{2,n})= 4^n r^{-2n} +\lambda^{2} (G_{2,n}) + 4^n \pi r^{-2n}$, a straightforward and sufficiently accurate upper bound is obtained by considering the convex hull of each of the $4^n$ connected components of $C^r$. This yields $\lambda^{2} (G_{2,n})<4^n \cdot 4 \cdot r^{-n} \cdot r^{-n}$. Thus, for all $n\in \NN$  it follows that

\begin{equation}
\begin{split}
\label{eqn:H2_vom_ersten_Satz}
\frac{\lambda^{2} (C^r_{\varepsilon_{2,n}})}{\varepsilon_{2,n}^{2-D_r}} &< \frac{4^n r^{-2n} \Big(1+4+\pi\Big)}{\Big(r^{-n}\Big)^2\Big(r^{-n}\Big)^{-\log_r(4)}} \\\ \\\
&= 5+\pi. 
\end{split}
\end{equation}

From Equation \eqref{eqn:H1_vom_ersten_Satz}, there exists an accumulation point $H_1^r$ of $\Big(\frac{\lambda^{2} (C^r_{\varepsilon_{1,n}})}{\varepsilon_{1,n}^{2-D_r}} \Big)_{n\in \NN}$ satisfying  $H^r_1 > \frac{1 + 3\sqrt{\frac{1}{2}} + \frac{\pi}{2}}{\frac{1}{2}\sqrt{2}^{\log_r(4)}}.$ Similarly, from Equation \eqref{eqn:H2_vom_ersten_Satz}, there exists an accumulation point $H^r_2$ of  $\Big(\frac{\lambda^{2} (C^r_{\varepsilon_{2,n}})}{\varepsilon_{2,n}^{2-D_r}} \Big)_{n\in \NN}$ satisfying $H^r_2< 5+\pi$. 

Noting  $\frac{1 + 3\sqrt{\frac{1}{2}} + \frac{\pi}{2}}{\frac{1}{2}\sqrt{2}^{\log_r(4)}}$ is as function for $r>1$ strictly increasing. Solving
 $$ \frac{1 + 3\sqrt{\frac{1}{2}} + \frac{\pi}{2}}{\frac{1}{2}\sqrt{2}^{\log_r(4)}} = 5+\pi \quad \Leftrightarrow \quad \frac{\log\Big(\frac{2+6\sqrt{\frac{1}{2}}+\pi}{5+\pi}\Big)}{\log(\sqrt{2})}=\log_r(4) \quad \Rightarrow \quad r \approx 29,4$$
leads to the conclusion that $H_1^r \gee H_2^r$ for $r\geq 30$, which implies  that $C^r$ is not Minkowski measurable for $r\geq 30$.
     
 \end{proof}

 \begin{remark}
    The sequences $\varepsilon_{1,n}=\sqrt{\frac{1}{2}}r^{-n}$ and $\varepsilon_{2,n}=r^{-n}$ do not satisfy the Equation \eqref{eqn:geeignete_Folgenwahl} for $r \in (2, 2+\sqrt{\frac{1}{2}})$.  Consequently, they are unsuitable for demonstrating the non-Minkowski measurability of $C^r$ for all $r>2$ using the approach from the proof above.
 \end{remark}

 \subsection{The case when $r < 30$}
 \label{ssec: case_r<30}
 
Considering the proof of Theorem \ref{Thm:case_r_grater_30}, to obtain an analogous result for all $ r > 2$, the sequences must be adjusted, and the estimates for the green part must be refined.  The following proof idea incorporates these adjustments, leading to an inequality that, although not rigorously proven,  has been numerically investigated without encountering any contradiction.

 \begin{conjecture}
     Let  $ r > 2$. Then $C^r$ is not Minkowski measurable 
 \end{conjecture}

\begin{proof}
Let $r > 2$ and consider two null sequences defined by $$\varepsilon_{1,n}:=\frac{1}{2}\frac{r-2}{r}\left(\frac{1}{r}\right)^{n-1} \qquad \text{and} \qquad \varepsilon_{2,n}:= \sqrt{\frac{1}{2}}\frac{r-2}{r}\left(\frac{1}{r}\right)^n.$$ For $r > 2$ $\varepsilon_{1,n}$ and $\varepsilon_{2,n}$ satisfy Equation \eqref{eqn:geeignete_Folgenwahl} by construction. Hence, $\lambda^{2} (C^r_{\varepsilon_{1,n}})$ is computed as follows:

\begin{equation}
\label{eqn:Formel_eins_Aus_Fall_r_grosser_zwei}
\begin{split}
\lambda^{2} (C^r_{\varepsilon_{1,n}}) &= \lambda^{2} (B_{1,n} \cup G_{1,n} \cup R_{1,n}) \\\ \\\
&= 4^n r^{-2n} + \lambda^{2} ( G_{1,n}) + 4^n \left( \frac{r-2}{r} \right)^2 \frac{r^2 \pi}{4} r^{-2n}.
\end{split}
\end{equation}

A sufficiently good lower bound for  $\lambda^2(G_{1,n})$ is shown in Figure \ref{fig:lower_bound_case_r_greater_two}.

 \begin{figure}[th]
	\centering
	\begin{tabular}{m{0.03\textwidth}m{0.87\textwidth}} 
	 $G_{1,n}$  & 
 
	\begin{tikzpicture}[rotate=90, transform shape,scale=0.8]
\filldraw[color=green!50, fill=green!40] (1,0) -- (3,0) -- (3,1.66) arc[start angle=0, end angle=40, radius= 1] coordinate(C) arc[start angle=-40, end angle=0, radius= 1] coordinate(D)-- (3,5) arc[start angle=0, end angle=25, radius= 6] coordinate(E) arc[start angle=-25, end angle=0, radius= 6] coordinate(F)--(3,12.07) arc[start angle=0, end angle=40, radius= 1] coordinate(G) arc[start angle=-40, end angle=0, radius= 1] coordinate(H)-- (3,15) -- (1,15);

\draw[dotted] (2.4,0) --(2.4,15);
\draw[dotted] (2.4,10.4) --(3,10.4);
\draw[dotted] (2.4,4.6) --(3,4.6);
\draw[dotted] (3,4.6) --(2.4,0);
\draw[dotted] (2.4,4.6) --(3,0);
\draw[dotted] (3,4.6) --(2.4,7.5);
\draw[dotted] (2.4,7.5) --(3,10.4);
\draw[dotted] (3,10.4) --(2.4,15);
\draw[dotted] (2.4,10.4) --(3,15);

	\end{tikzpicture} \\

 \phantom{x}&
\\
	$G'_{1,n}$ &
 
	\begin{tikzpicture}[rotate=90, transform shape,scale=0.8]
\filldraw[color=green!50, fill=green!40] (1,0)--(3,0)--(2.4,4.6)--(3,4.6)--(2.4,7.5)--(3,10.4)--(2.4,15)--(1,15)--(1,0);
\filldraw[color=green!50, fill=green!40] (1,0)--(2.4,0)--(3,4.6)--(2.4,7.5)--(3,10.4)--(2.4,10.4)--(3,15)--(1,15)--(1,0);

\draw[dotted] (2.4,0) --(2.4,15);
\draw[dotted] (2.4,10.4) --(3,10.4);
\draw[dotted] (2.4,4.6) --(3,4.6);
\draw[dotted] (3,4.6) --(2.4,0);
\draw[dotted] (2.4,4.6) --(3,0);
\draw[dotted] (3,4.6) --(2.4,7.5);
\draw[dotted] (2.4,7.5) --(3,10.4);
\draw[dotted] (3,10.4) --(2.4,15);
\draw[dotted] (2.4,10.4) --(3,15);
	\end{tikzpicture}

	\end{tabular}
	\caption{Representation of the lower bound $G'_{1,n}$ of $G_{1,n}$.}
	\label{fig:lower_bound_case_r_greater_two}
\end{figure}
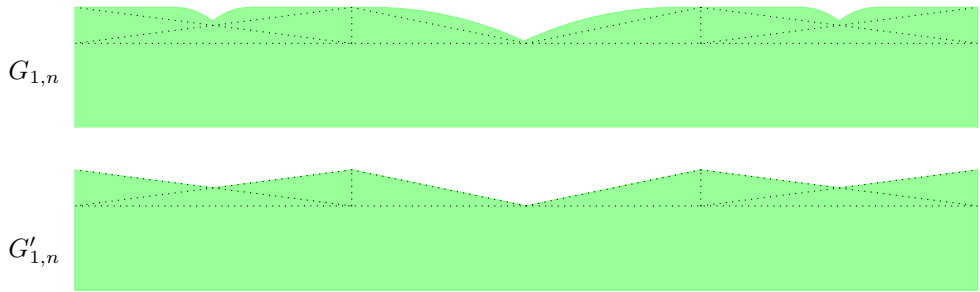  

Figure \ref{fig:lower_bound_for_r_greater_30} shows that the estimation depicted in Figure \ref{fig:lower_bound_case_r_greater_two} is indeed a lower bound for all $r>2$ (see also the paragraph directly before Equation \eqref{eqn:H1_vom_ersten_Satz}). Expressed in formulas, first calculate the smallest distance $h$ from  $\partial C^r_{\varepsilon_{1,n}}$ to $B_{1,n}$, see Figure \ref{fig:berechnung_von_h}. By the Pythagorean theorem, $h$ is given by

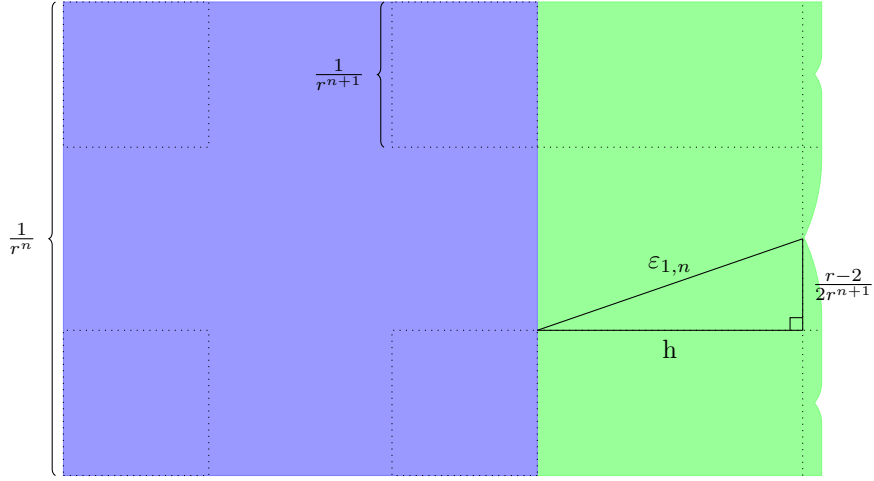
\begin{figure}[th]
	\centering

	\begin{tikzpicture}[scale=0.42]

\filldraw[color=green!50, fill=green!40] (0,0) -- (9,0) -- (9,1.66) arc[start angle=0, end angle=40, radius= 1] coordinate(C) arc[start angle=-40, end angle=0, radius= 1] coordinate(D)-- (9,5) arc[start angle=0, end angle=25, radius= 6] coordinate(E) arc[start angle=-25, end angle=0, radius= 6] coordinate(F)--(9,12.07) arc[start angle=0, end angle=40, radius= 1] coordinate(G) arc[start angle=-40, end angle=0, radius= 1] coordinate(H)-- (9,15) -- (0,15);

\filldraw[color=blue!50, fill=blue!40] (0,0) -- (-15,0)--(-15,15)--(0,15)--(0,0);
\draw[dotted] (0,0) -- (-4.6,0)--(-4.6,4.6)--(0,4.6)--(0,0);
\draw[dotted] (-15,0) -- (-15+4.6,0)--(-15+4.6,4.6)--(-15,4.6)--(-15,0);
\draw[dotted] (-15,15) -- (-15+4.6,15)--(-15+4.6,15-4.6)--(-15,15-4.6)--(-15,15);
\draw[dotted] (-4.6,15-4.6)--(0,15-4.6)--(0,15) -- (-4.6,15);
\draw[dotted]  (-4.6,15)--(-4.6,15-4.6);

\draw[decoration={brace,mirror,raise=3pt},decorate]
  (-4.6,15) -- node[midway,left] {$\frac{1}{r^{n+1}}\ \ $} (-4.6,15-4.6);

  \draw[decoration={brace, mirror,raise=3pt},decorate]
  (-15,15) -- node[midway,left] {$\ \ \frac{1}{r^{n}}\ \ $} (-15,0);

\draw[dotted] (0,4.6) -- (9,4.6);
\draw[dotted] (0,15-4.6) -- (9,15-4.6);
\draw[dotted] (8.4,0) -- (8.4,15);
\draw (8,4.6) -- (8, 5) -- (8.4, 5);
\draw (8.4,7.5) -- (0,4.6) node[midway,above] {$\varepsilon_{1,n}$};
\draw (0,4.6) -- (8.4, 4.6) node[midway,below] {h};
\draw (8.4, 4.6) -- (8.4, 7.5) node[midway,right] {$\frac{r-2}{2r^{n+1}}$};

	\end{tikzpicture}
	\caption{The smallest distance $h$ from $\partial C^r_{\varepsilon_{1,n}}$ to $B_{1,n}$.}
	\label{fig:berechnung_von_h}
\end{figure}

\begin{equation}
\begin{split}
h &= \sqrt{\varepsilon_{1,n}^2 - \Bigg( \frac{r-2}{2r^{n+1}} \Bigg )^2} \\\ \\\
&=\frac{1}{2}\frac{r-2}{r}\left(\frac{1}{r}\right)^n\sqrt{r^2-1}.
\end{split}
\end{equation}

Therefore,
\begin{equation*}
\begin{split}
\lambda^2(G_{1,n}) &> 4^n \cdot 4 \Bigg [ r^{-n} \cdot h + \frac{r-2}{2 r^{n+1}} \cdot (\varepsilon_{1,n} - h) + 2\cdot \frac{3}{4r^{n+1}} \cdot (\varepsilon_{1,n} - h) \Bigg ] \\\ \\\
&= 4^{n+1} \cdot \frac{1}{4}r^{-2n}\cdot \Bigg [ 2\frac{r-2}{r}\sqrt{r^2-1}+3\frac{r-2}{r}-\frac{3}{r}\frac{r-2}{r}\sqrt{r^2-1}+\left(\frac{r-2}{r}\right)^2\left(r-\sqrt{r^2-1}\right) \Bigg ].
\end{split}
\end{equation*}

Equation (\ref{eqn:Formel_eins_Aus_Fall_r_grosser_zwei}) gives, for all $n \in \NN$,

\begin{equation}
\begin{split}
\label{eqn:H1_vom_zweiten_Satz}
\frac{\lambda^{2} (C^r_{\varepsilon_{1,n}})}{\varepsilon_{1,n}^{2-D_r}} &= \frac{4^n r^{-2n} + \lambda^{2} ( G_{1,n}) + 4^n \left( \frac{r-2}{r} \right)^2 \frac{r^2 \pi}{4} r^{-2n}}{\varepsilon_{1,n}^{2-D_r}} \\\ \\\
&>\frac{1+2\frac{r-2}{r}\sqrt{r^2-1}+3\frac{r-2}{r}-\frac{3}{r}\frac{r-2}{r}\sqrt{r^2-1}+\left(\frac{r-2}{r}\right)^2\left(r-\sqrt{r^2-1}+\frac{r^2\pi}{4}\right)}
{\frac{1}{4}\left(r-2\right)^2\left(\frac{1}{2}\right)^{-\log_r 4}\left(r-2\right)^{-\log_r 4}} \\\ \\\
&:= f_1(r).
\end{split}
\end{equation}

For $\lambda^{2} (C^r_{\varepsilon_{2,n}})$ the following is obtained:
\begin{equation}
\label{eqn:Formel_zwei_Aus_Fall_r_grosser_zwei}
\begin{split}
\lambda^{2} (C^r_{\varepsilon_{2,n}}) &= \lambda^{2} (B_{2,n} \cup G_{2,n} \cup R_{2,n}) \\\ \\\
&= 4^n r^{-2n} +\lambda^{2} (G_{2,n}) + 4^n \frac{\pi}{2}\left( \frac{r-2}{r} \right)^2 r^{-2n}.
\end{split}
\end{equation}

An upper bound for $\lambda^2 (G_{2,n})$ is given by the convex hull of each of the $4^n$ connected components of $C^r_{\varepsilon_{2,n}}$:

$$\lambda^{2} (G_{2,n})<4^n \cdot 4 \cdot \sqrt{\frac{1}{2}}\frac{r-2}{r} r^{-n} \cdot r^{-n}= 4^n \sqrt{8} \frac{r-2}{r}r^{-2n}.$$

Together with Equation \eqref{eqn:Formel_zwei_Aus_Fall_r_grosser_zwei}, the following holds for all $n \in \NN$:

\begin{equation}
\begin{split}
\label{eqn:H2_vom_zweiten_Satz}
\frac{\lambda^{2} (C^r_{\varepsilon_{2,n}})}{\varepsilon_{2,n}^{2-D_r}} &= \frac{4^n r^{-2n} + \lambda^{2} ( G_{2,n}) +  4^n \frac{\pi}{2}\left( \frac{r-2}{r} \right)^2 r^{-2n}}{\varepsilon_{2,n}^{2-D_r}} \\\ \\\
&< \frac{1+\sqrt{8}\frac{r-2}{r}+\frac{\pi}{2}\left(\frac{r-2}{r}\right)^2}{\frac{1}{2}\left(\frac{r-2}{r}\right)^2\sqrt{\frac{1}{2}}^{-\log_r 4}\left(\frac{r-2}{r}\right)^{-\log_r 4}} \\\ \\\
&:= f_2(r).
\end{split}
\end{equation}

From Equation \eqref{eqn:H1_vom_zweiten_Satz} there exists an accumulation point $H^r_1$ of $\Big(\frac{\lambda^{2} (C^r_{\varepsilon_{1,n}})}{\varepsilon_{1,n}^{2-D_r}} \Big)_{n\in \NN}$ such that $H^r_1 >f_1(r)$. Similarly, from Equation \eqref{eqn:H2_vom_zweiten_Satz} there exists an accumulation point $H^r_2$ of $\Big(\frac{\lambda^{2} (C^r_{\varepsilon_{2,n}})}{\varepsilon_{2,n}^{2-D_r}} \Big)_{n\in \NN}$ satisfying $H^r_2 < f_2(r)$.

Hence, if $f_1(r) \geq f_2(r)$, then
\begin{equation}
\label{eqn:f_1_groesser_f_2}
H^r_1>f_1(r) \geq f_2(r) > H^r_2
\end{equation}
which implies that $C^r$ is not Minkowsky measurable for $r > 2$. \\

Analytically, the inequality $f_1(r)\geq f_2(r)$ has not yet been confirmed. However, a Python code - using a step size of $10^{-4}$ -  computes Equation \eqref{eqn:f_1_groesser_f_2} for $r\in [2.0001, 30]$ and verifies it, see Figure \ref{fig:Plots_of_funktions}. Note that in Theorem \ref{Thm:case_r_grater_30} the non-Minkowski measurability of $C^r$ for  $r \geq 30$ has already been established. \\

 \begin{figure}
     \centering
     \centering
	\begin{tabular}{ c @{\qquad \qquad} c }
         \includegraphics[width=0.42\linewidth]{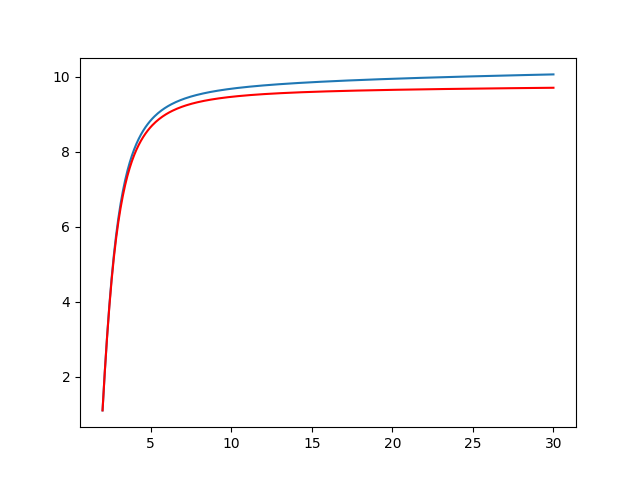}
        &
        \includegraphics[width=0.42\linewidth]{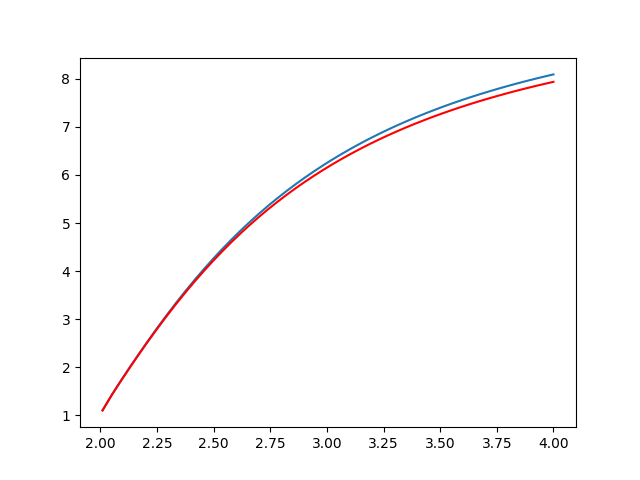}

	\end{tabular}
    
     \caption{Plots of the functions $f_1$ in blue and $f_2$ in red}
     \label{fig:Plots_of_funktions}
 \end{figure}

\end{proof}

\begin{remark}
The case $r=3$ is known from \emph{\cite[ch.~5.1]{Beck}}, where Matthias Beck computed several examples in support of the Lapidus conjecture.   
\end{remark}

\section*{Acknowledgements} We would like to express our gratitude to Nico Heinzmann and Steffen Winter for their careful proofreading. Their insightful comments greatly improved the clarity and accuracy of this work.  In particular, thanks are due to Steffen Winter for simplifying the proof of Proposition  \ref{prop:C^r_pluri}. Additionally, we acknowledge the support of the European Social Fund (ESF) and TU Chemnitz.

\bibliography{refs}

\end{document}